\documentclass[12pt]{article}

\usepackage{amsmath, amsfonts,epsfig}
\usepackage{url}

\newtheorem{theorem}{Theorem}
\newtheorem{lemma}{Lemma}
\newtheorem{remark}{Remark}

\newtheorem{corollary}{Corollary}
\newtheorem{proposition}{Proposition}

\newcommand\pf{{\bf Proof.} }
\newcommand\qed{\hfill$\Box$}

\newcommand\Gm{\Gamma}
\newcommand\Om{\Omega}

\newcommand\al{\alpha}
\newcommand\be{\beta}
\newcommand\gm{\gamma}
\newcommand\dl{\delta}
\newcommand\ep{\varepsilon}

\title{Graphs $4_n$ that are isometrically embeddable in hypercubes}

\author{Michel DEZA\\
        \normalsize   LIGA, ENS, Paris and Institute of Statistical 
                      Mathematics, Tokyo\\
\and
        Mathieu DUTOUR-SIKIRIC\thanks{Research financed by EC's IHRP 
                              Programme, within the Research Training 
                              Network ``Algebraic Combinatorics in 
                              Europe,'' grant HPRN-CT-2001-00272.}\\
        \normalsize  LIGA, ENS, Paris and Hebrew University, 
                     Jerusalem\\
\and
        Sergey SHPECTOROV\thanks{Research partly supported by an 
                                 NSA grant.}\\
        \normalsize Bowling Green State University, Bowling Green\\
}

\begin{document}

\maketitle

\begin{abstract}
A connected $3$-valent plane graph, whose faces are $q$- or $6$-gons
only, is called a {\em graph $q_n$}. We classify all graphs
$4_n$, which are isometric subgraphs of a $m$-hypercube $H_m$.
\end{abstract}

\section{Introduction}

A $3$-valent $n$-vertex plane graph is denoted $q_n$
if it has only $q$- and $6$-gonal faces.
The graphs $5_n$ correspond to {\em fullerenes}, well-known in 
Organic Chemistry; graphs $4_n, 3_n$ are also used there.

See \cite{DG2} for $\ell_1$-embedding of graphs $q_n$.
Denote by $Aut$ the automorphism group of given graph $q_n$.

The (vertex-set of) hypercube $H_m=\{0,1\}^m$ is a metric space under 
the distance
$d(x,y)=\sum_i |x_i-y_i|$. 

%In this paper we identify those graphs $4_n$, which
%are isometric subgraphs of a hypercube $H_m$.

See \cite{DGS} (especially, Chapter 13) and \cite{zig2} on 
graphs $q_n$.

A scale $\lambda$
embedding $\phi$ of a graph $G$ into a hypercube $H_m$ is a
mapping $G\mapsto H_m$, such that
\begin{equation*}
\lambda d_G(x,y)=d(\phi(x), \phi(y))
\end{equation*}
with $d_G$ being the path-distance on $G$.
It was shown in \cite{AsDe} that $d_G$ (moreover, any finite 
rationally-valued metric space) is {\em $l_1$-embeddable}
(i.e. embeds isometrically into some $l_1^k$) if and only if it is scale 
$\lambda$ embeddable in $H_m$ for some $\lambda$ and $m$.

The cases $\lambda=1$ or $2$ mean exactly that $G$ is an isometric
subgraph of some hypercube $H_m$ or, respectively,
of some {\em half-cube} $\frac{1}{2}H_m$, where
\begin{equation*}
\frac{1}{2}H_m=\{x \in \{0,1 \} ^m : \sum_i x_i\mbox{~~is~even} \}
\end{equation*}
and two vertices are adjacent if and only if $d(x,y)=2$.

Any $l_1$-embeddable graph satisfy (see \cite{D}, \cite{DL}) the 
following {\em $5$-gonal inequality}:
\begin{equation*}
\begin{array}{rcl}
d(a,b)+d(x,y)+d(x,z)+d(y,z)&\leq& d(a,x)+d(a,y)+d(a,z)\\
                           &+   &d(b,x)+d(b,y)+d(b,z).
\end{array}
\end{equation*}

For a {\em bipartite} graph (so, in particular, for any $4_n$)
the following three conditions are equivalent: 
$l_1$-embeddability, being an isometric 
subgraph of a hypercube and satisfy all $5$-gonal inequalities (see 
\cite{Av} and Chapter 19.2 in \cite{DL}).

\section{Main Theorem}

\begin{theorem} \label{main}
If a $4_n$ graph $\Gm$ is isometrically embeddable in a hypercube then
$\Gm$ is one of the following graphs: the Cube ($O_h$), the $6$-gonal prism ($D_{6h}$), the truncated Octahedron ($O_h$), the chamfered Cube ($O_h$), or the twisted chamfered Cube ($D_{3h}$).
\end{theorem}

The pictures of last three polyhedra are given below in Lemmata 
\ref{x_and_y_equal}, \ref{x_and_y_different} and \ref{star_star},
respectively.

We prove this theorem in a sequence of lemmas. Suppose $\phi$ is an
isometric embedding of $\Gm$ into a hypercube $H_m$. We realize $H_m$ as
the graph whose vertices are all subset of a base set $\Om$, $|\Om|=m$. The
elements of $\Om$ will be called the {\em coordinates} of the hypercube
$H_m$. Two vertices-subsets $A$ and $B$ are {\em adjacent} if and only if
$A$ and $B$ differ in just one coordinate, that is, if the {\em symmetric
difference} $A\triangle B=(A\setminus B)\cup(B\setminus A)$ is of size
one. Because of this, every edge of $H_m$ naturally carries a label, which
is one of the coordinates. More precisely, the edge $AB$ is labelled with
the unique element from $A\triangle B$.

In general, if $A$ and $B$ are arbitrary vertices of $H_m$ then the
distance between $A$ and $B$ is exactly $|A\triangle B|$. Furthermore, the
labels on the consecutive edges along an arbitrary shortest path from $A$
to $B$ are simply the elements of $A\triangle B$, appearing in some order. 

If we now consider an arbitrary (that is, not necessarily shortest) path
from $A$ to $B$ then the labels on the edges along this path may or may not
belong to $A\triangle B$, and they may or may not repeat. However, if we
count the labels modulo two, that is, if we only take the labels that
appear {\em oddly-many} times along the path, then they are exactly the
elements of $A\triangle B$.

Looking at this from a different angle, recall that a path is called {\em
geodesic} if it is one of the shortest paths from one of its ends to the
other. It follows from the above that if a path in $H_m$ is geodesic then
the labels along the path do not repeat, and the (unordered) set of lebels
along the path, which we will refer to as the {\em type} of the geodesic
path, can be found as the symmetric difference $A\triangle B$, where $A$
and $B$ are the ends of the path.

Now we turn to the graph $\Gm$ and its isometric embedding $\phi$. Since
$\phi$ takes edges of $\Gm$ to edges of $H_m$, every edge of $\Gm$ gets a
{\em label}, which is an element of $\Om$.  Similarly, since $\phi$ takes
geodesic paths of $\Gm$ to geodesic paths of $H_m$, every geodesic path of
$\Gm$ gets a {\em type}, which coincides with the type of its image under
$\phi$, and which is, simply, the set of labels on the edges along the
path. The type of a geodesic path from a vertex $\al$ to a vertex $\be$ of
$\Gm$ can be computed as $\phi(\al)\triangle\phi(\be)$. If $\phi(\al)$ and
$\phi(\be)$ are not known, but we know the labels along some (not
necessarily geodesic) path from $\al$ to $\be$, then we can count the
labels on this path modulo two. The resulting set will be the type of every
geodesic path connecting $\al$ and $\be$. This simple trick will be an
important tool in what follows.

Lastly, we will use small letters for the labels---elements of $\Om$. We
will use (as above) small greek letters for the vertices of $\Gm$.

We now start the proof.

\begin{lemma} \label{claw}
The labels on adjacent edges of $\Gm$ are never equal.
\end{lemma}

\pf Indeed, $\Gm$ is bipartite, since $H_m$ is. This means that any two
adjacent edges of $\Gm$ form a geodesic path.\qed

\medskip
This means that we have the following picture around every vertex.
\begin{center}
\epsfig{file=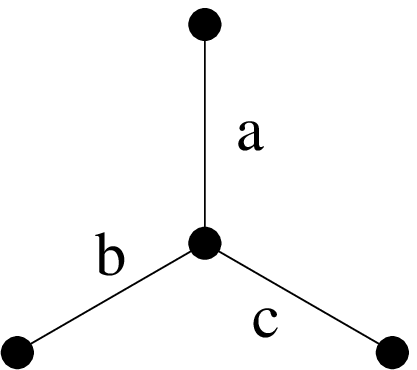,width=2cm}
\end{center}
Here the labels $a$, $b$, and $c$ are pairwise distinct.

\begin{lemma} \label{four-cycles are isometric}
Every $4$-cycle in $\Gm$ is isometric.
\end{lemma}

\pf Since $\Gm$ is bipartite, the distance between the opposite vertices on
the $4$-cycle cannot be one.\qed

\begin{corollary} \label{labeling of four-cycles}
The edges of every $4$-cycle carry labels as follows:
\begin{center}
\epsfig{file=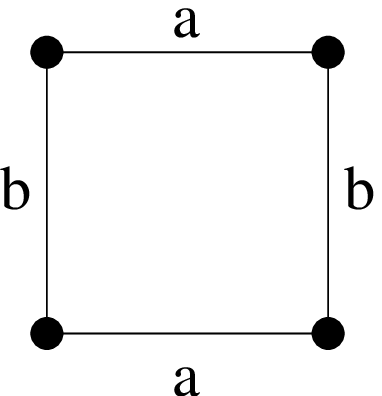,width=2cm}
\end{center}
The labels $a$ and $b$ are distinct.
\end{corollary}

\pf Let $\al$ and $\be$ be nonadjacent vertices of the $4$-cycle.  Within
the cycle, we have two different geodesic paths from $\al$ to $\be$.  Both
paths must have the same type, say, $\{a,b\}$. Finally, by Lemma
\ref{claw}, if one of the paths starts with the label $a$ then the other
one must start with $b$, and vise versa.\qed

\begin{lemma} \label{four-cycles are faces}
Every $4$-cycle in $\Gm$ is a face.
\end{lemma}

\pf Suppose a $4$-cycle is not a face. Then among the edges adjacent to the
cycle there are both edges pointing inside the cycle and outside the
cycle. It is easy to see that we have one of the following two situations:
\begin{center}
\epsfig{file=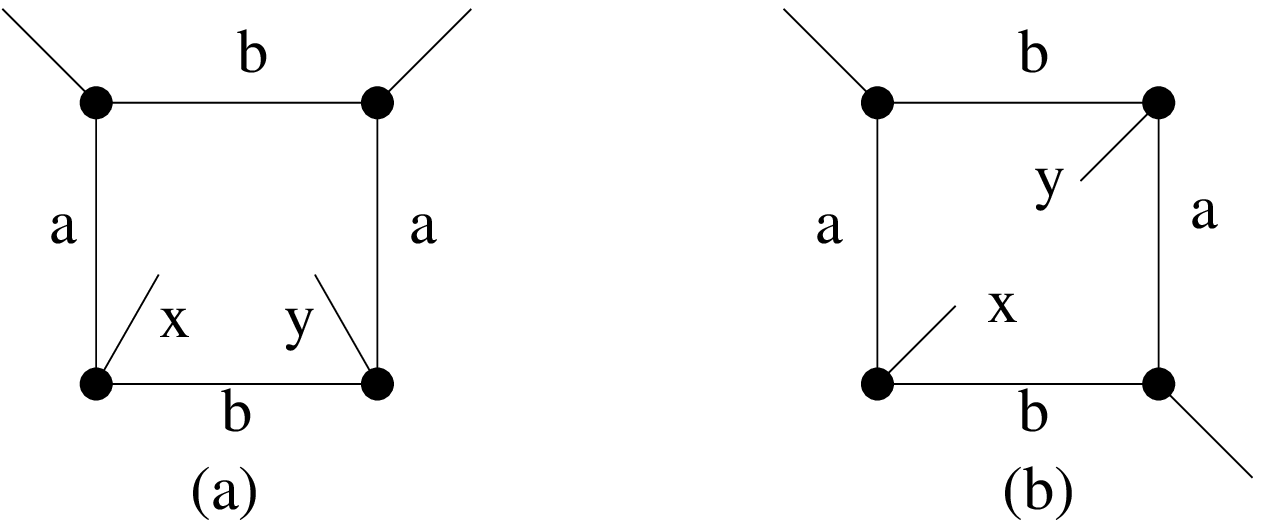,width=6cm}
\end{center}

\medskip
{\em Case 1.} Here the consecutive edges $x$, $a$, $b$, $a$, and $y$ lie on
a face. This face is $6$-gonal, and we end up with the following picture.
\begin{center}
\epsfig{file=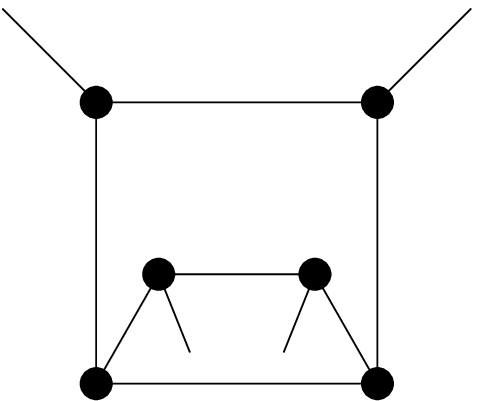,width=3cm}
\end{center}
Notice that the original configuration recurs inside the big $4$-gon, so we
get a contradiction with the finiteness of $\Gm$.

\medskip
{\em Case 2:} Suppose we have the situation as in picture (b). The
consecutive edges labeled $x$, $a$, $b$, and $y$ lie on the same
face. Clearly, $x,y\not\in\{a,b\}$. So this face is $6$-gonal. Similarly,
for the consecutive edges $x$, $b$, $a$, and $y$. Thus, we have the
following picture.
\begin{center}
\epsfig{file=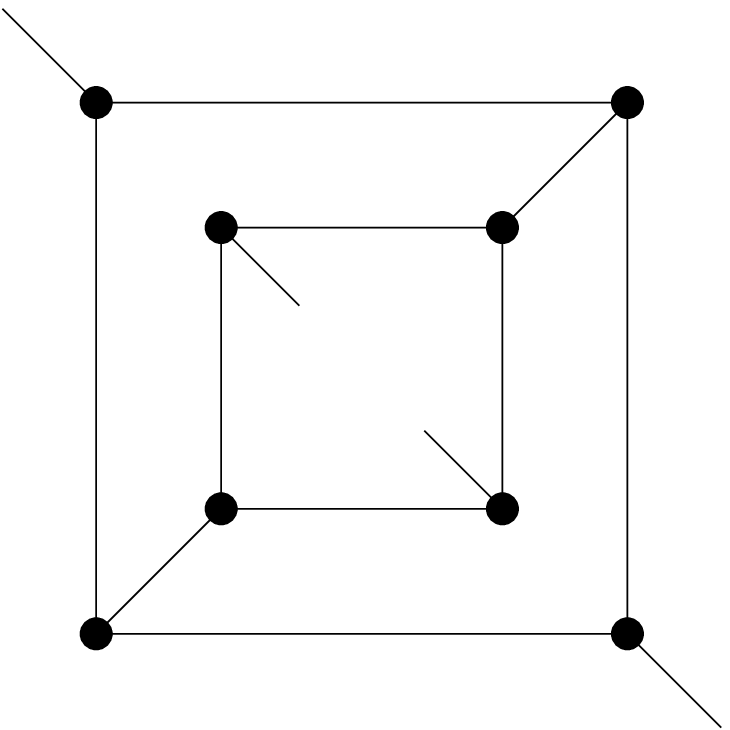,width=3cm}
\end{center}
Again, the situation recurs inside the original $4$-gon. Since $\Gm$ is
finite, we get a contradiction.\qed

\begin{lemma} \label{six-gonal faces are isometric}
Every $6$-gonal face is isometric.
\end{lemma}

\pf Suppose not. Then two opposite vertices of the face must be at distance
one. (Since $\Gm$ is bipartite, the distance cannot be two.) This leads to
two $4$-cycles, each of which must be a face by Lemma \ref{four-cycles are
faces}. We get a contradiction, since the original $6$-gonal face and each
of the $4$-gonal faces share a $3$-path.\qed

\begin{corollary} \label{labels on six-faces}
Every $6$-gonal face carries labels as follows.
\begin{center}
\epsfig{file=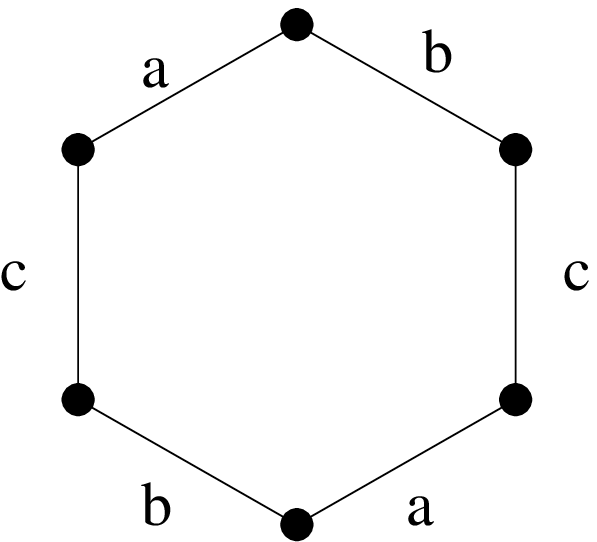,width=3cm}
\end{center}
Here the labels $a$, $b$, and $c$ are pairwise distinct.
\end{corollary}

\pf Taking two opposite vertices as $\al$ and $\be$, we see that the two
paths between $\al$ and $\be$ in this $6$-cycle must have the same type,
say $\{a,b,c\}$. On the other hand, every $3$-path in this $6$-cycle is
geodesic, which means that the labels on it must be distinct. Thus, the
opposite edges must have the same label.\qed

\medskip
From this point on we are looking at various subcases.

\begin{lemma} \label{three four-gons}
If the three faces meeting at a vertex are $4$-gons then $\Gm$ is the Cube
graph.
\end{lemma}

\pf We have the following picture.
\begin{center}
\epsfig{file=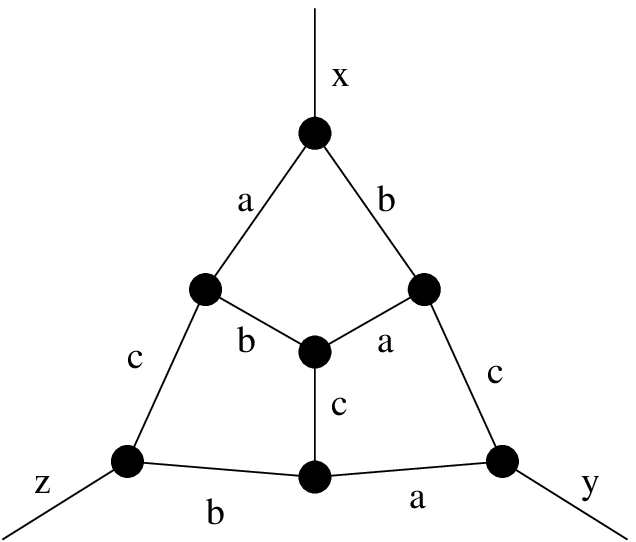,width=4cm}
\end{center}
Clearly, $x\ne a,b$. If $x=c$ then, since $x$, $b$, $c$, and $y$ label
consecutive edges of a face, we get that that face is a $4$-gon and hence
$y=b$. Using the same argument for the other two sides of the picture, we
obtain a Cube.

So now suppose that $x\ne c$, and similarly, $y\ne a,b,c$ and $z\ne
a,b,c$. Then on all three sides of the picture we get $6$-gonal faces,
yielding $x=y=z$. Now the following larger picture can be drawn, where the
additional labels are derived from Corollary \ref{labels on six-faces}.
\begin{center}
\epsfig{file=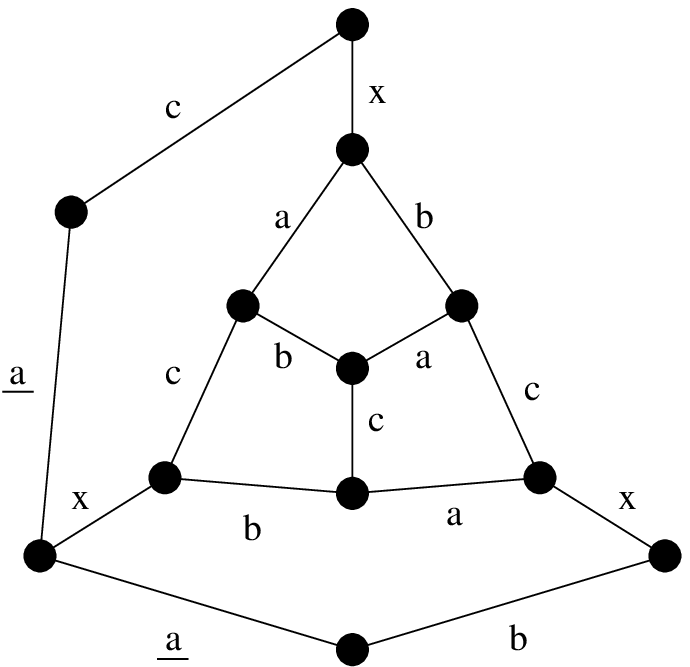,width=4cm}
\end{center}
However, two edges incident to the same vertex (bottom left) cannot have
the same label $a$, a contradiction.\qed

\begin{lemma} \label{four-gons next to each other}
If $\Gm$ contains two $4$-gonal faces next to each other then $\Gm$ is
either a Cube, or a $6$-gonal prism.
\end{lemma}

\pf We start with the following picture.
\begin{center}
\epsfig{file=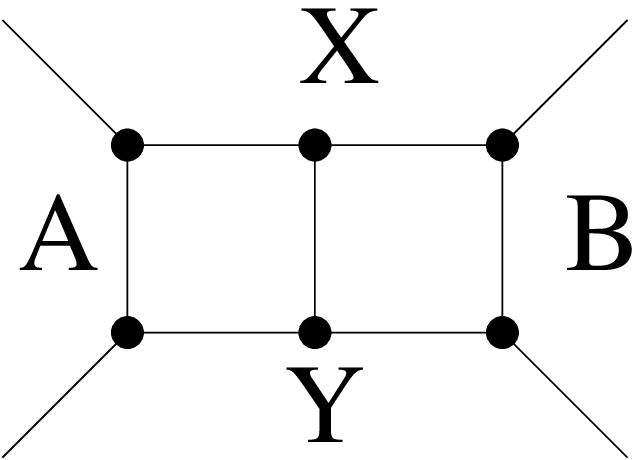,width=3cm}
\end{center}
If either of the faces $X$ and $Y$ is $4$-gonal then $\Gm$ is a Cube by
Lemma \ref{three four-gons}. So we will now assume that both $X$ and $Y$
are $6$-gonal. We claim that both $A$ and $B$ are $4$-gonal. By
contradiction, suppose that $A$ is $6$-gonal. Then, using Corollary
\ref{labels on six-faces}, we see that $B$ is also $6$-gonal and the
picture is as follows, where the labels $a$, $b$, $c$, $x$, and $y$ are
pairwise distinct and $u$ is an unknown label.
\begin{center}
\epsfig{file=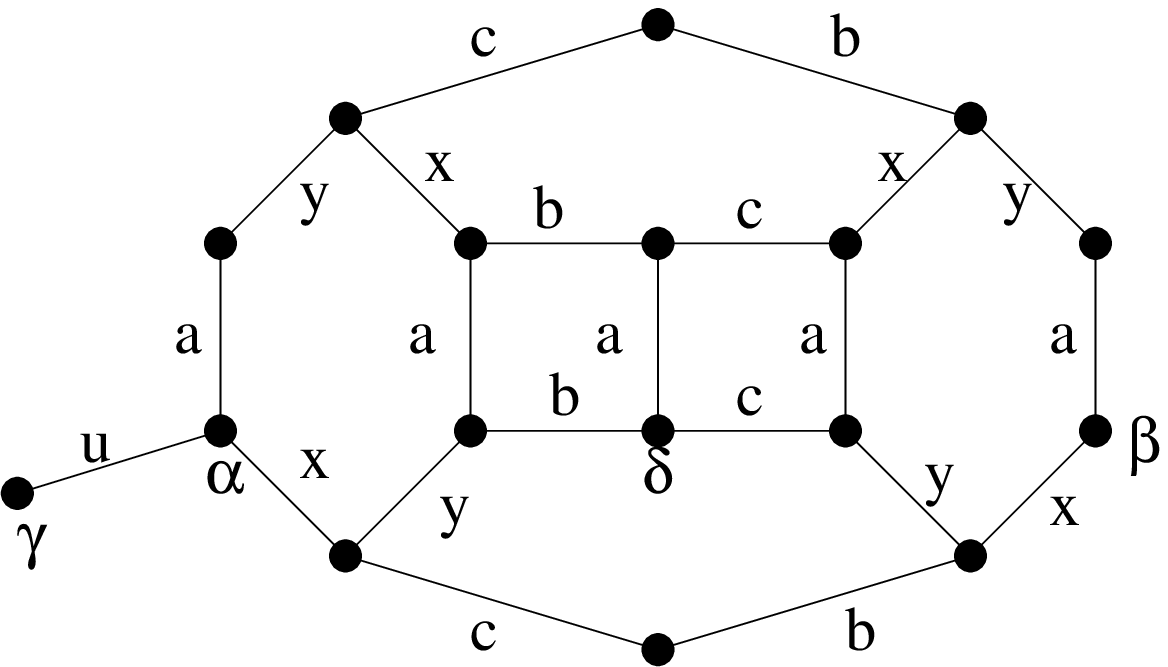,width=6cm}
\end{center}
We can now compute, counting modulo two the labels along an arbitrary path
from $\al$ to $\be$, that the type of every shortest path from $\al$ to
$\be$ is $\{b,c\}$. This means that $u=b$ or $u=c$. If $u=b$ then the type
of the shortest path from $\gm$ to $\dl$ is $\{x,y\}$, a contradiction
since $\dl$ has neither $x$, nor $y$ on the edges next to it. Thus,
$u=c$. Similarly, the third edge at $\be$ has label $b$. Furthermore, the
two new edges (at $\al$ and at $\be$) lead to the same new vertex, the
common neighbor of $\al$ and $\be$. The same reasoning can be applied to
the vertices above $\al$ and $\be$, which leads us to the following
picture.
\begin{center}
\epsfig{file=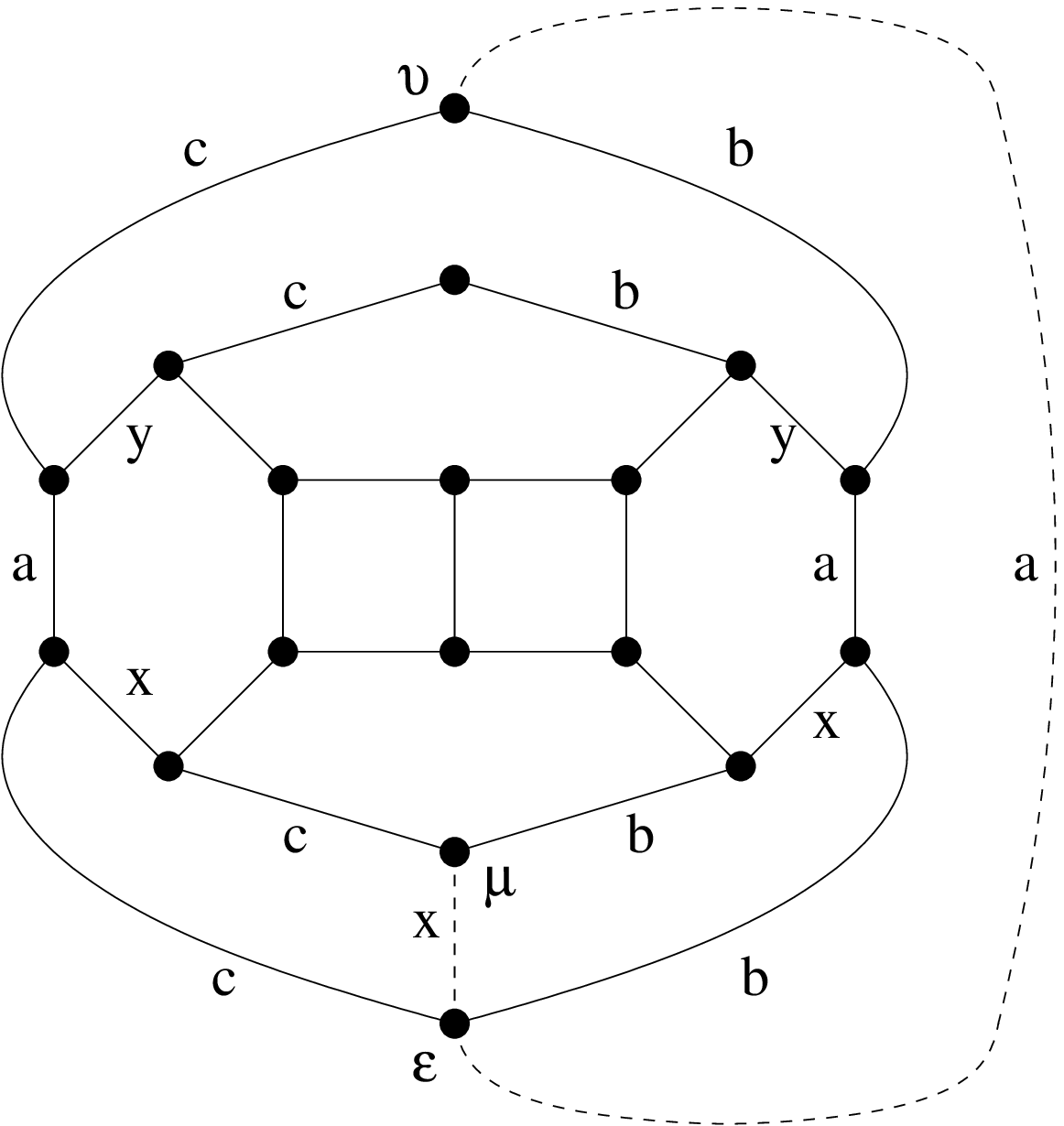,width=6cm}
\end{center}
Now, in this new picture the vertex $\ep$ must be adjacent to $\mu$ by an
edge labeled $x$, and to $\nu$ by an edge labeled $a$. This is a
contradiction since the valency of the graph is three. The contradiction
proves our claim.

Thus, returning to the first picture in this proof, we have that both faces
$A$ and $B$ are $4$-gons. Since this applies to any two adjacent $4$-gonal
faces, we end up with a prism.\qed

\medskip
From now on {\em we assume that $\Gm$ contains no adjacent $4$-gons}.  In
particular, every $4$-gonal face is surrounded only by $6$-gonal faces.

\begin{lemma} \label{six-cycles are faces}
Every $6$-cycle in $\Gm$ is a face.
\end{lemma}

\pf Consider a $6$-cycle in $\Gm$. Since $\Gm$ contains no adjacent
$4$-gons, the $6$-cycle is isometric. So the labels on this $6$-cycles are
as in Corollary \ref{labels on six-faces}.

Suppose first that the third edges at two consecutive vertices of the cycle
are on the same side of it. Then we have the following picture.
\begin{center}
\epsfig{file=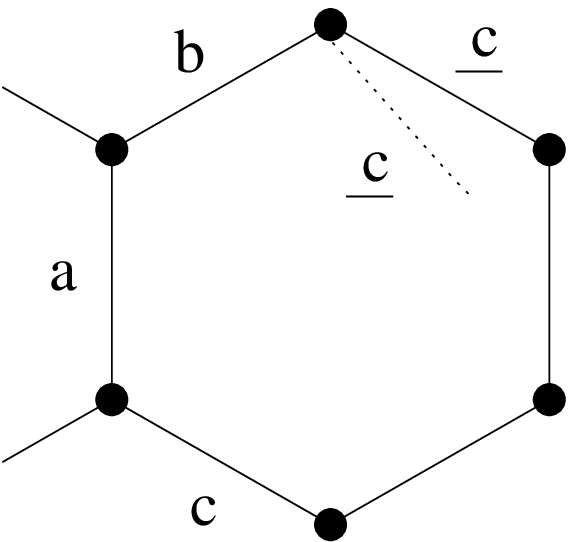,width=3cm}
\end{center}
The consecutive edges $c$, $a$, and $b$ (starting from the bottom center
vertex) belong to a face, which must be $6$-gonal. However, in that case,
the next edge on the face must be labelled $c$, which means that the
original cycle and the face share the next edge as well. Continuing in the
same manner, we obtain that our cycle is the face.

This leaves us with the case where the third edges along the cycle go to
alternate sides. So the picture is as follows.
\begin{center}
\epsfig{file=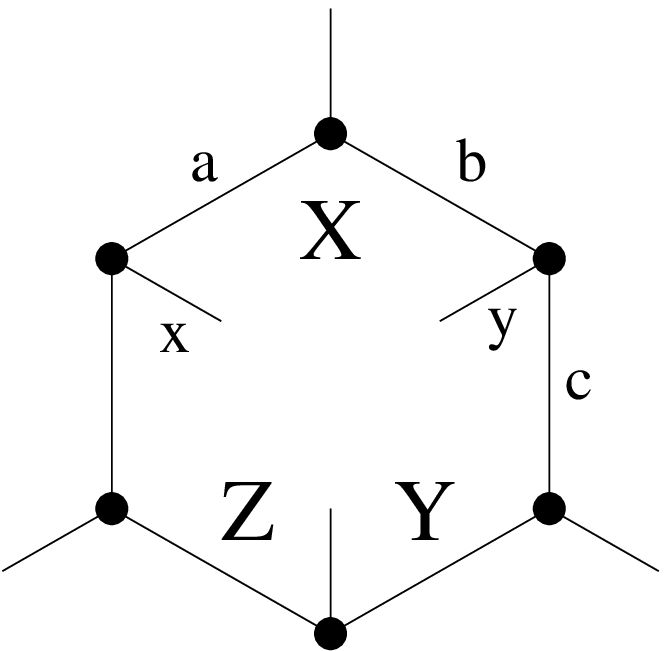,width=4cm}
\end{center}
If the face $X$ is $4$-gonal then $x=b$ and $y=a$. This immediately yields
that the faces $Y$ and $Z$ are also $4$-gonal and $\Gm$ is a Cube; a
contradiction to our assumption. Thus, $X$, $Y$, and $Z$ are $6$-gonal.
\begin{center}
\epsfig{file=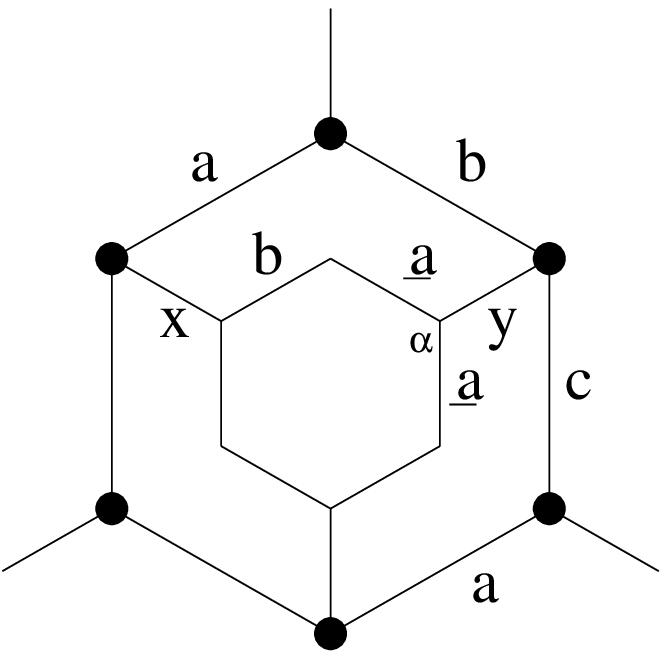,width=4cm}
\end{center}
However, in this case, two edges at the vertex $\al$ are both labelled with
$a$; a contradiction.\qed

\medskip
In the next two lemmas we analyze the case where $\Gm$ contains the
following configuration.
\begin{center}
\epsfig{file=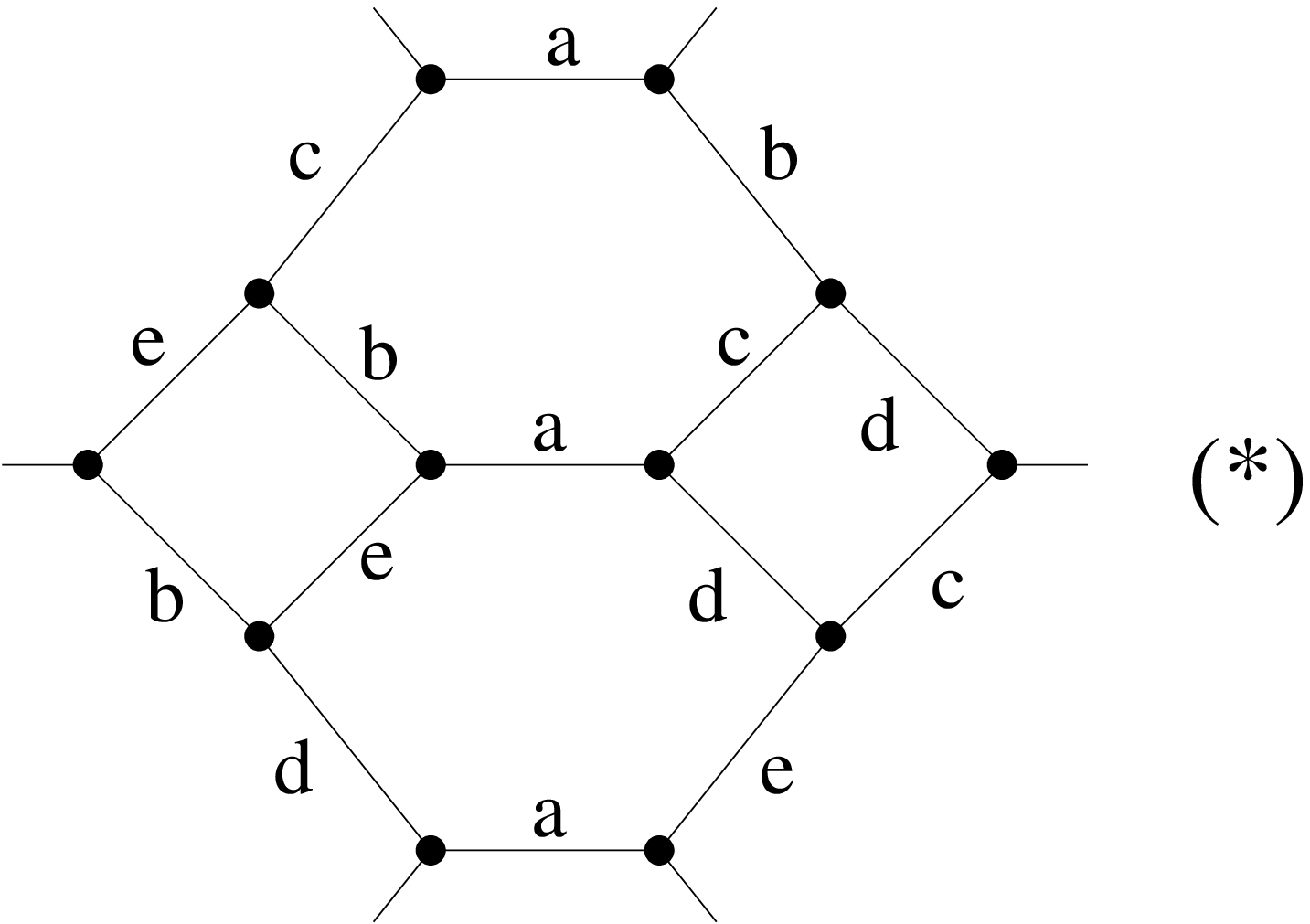,width=5cm}
\end{center}
Notice that the two faces in the center {\em must} be $6$-gonal due to our
assumption, since they border $4$-gonal faces. For the same reason, the
faces $X$, $Y$, $Z$, and $T$ below are also $6$-gonal, which gives us just
two new labels, $x$ and $y$. Clearly, $x,y\not\in\{a,b,c,d,e\}$.
\begin{center}
\epsfig{file=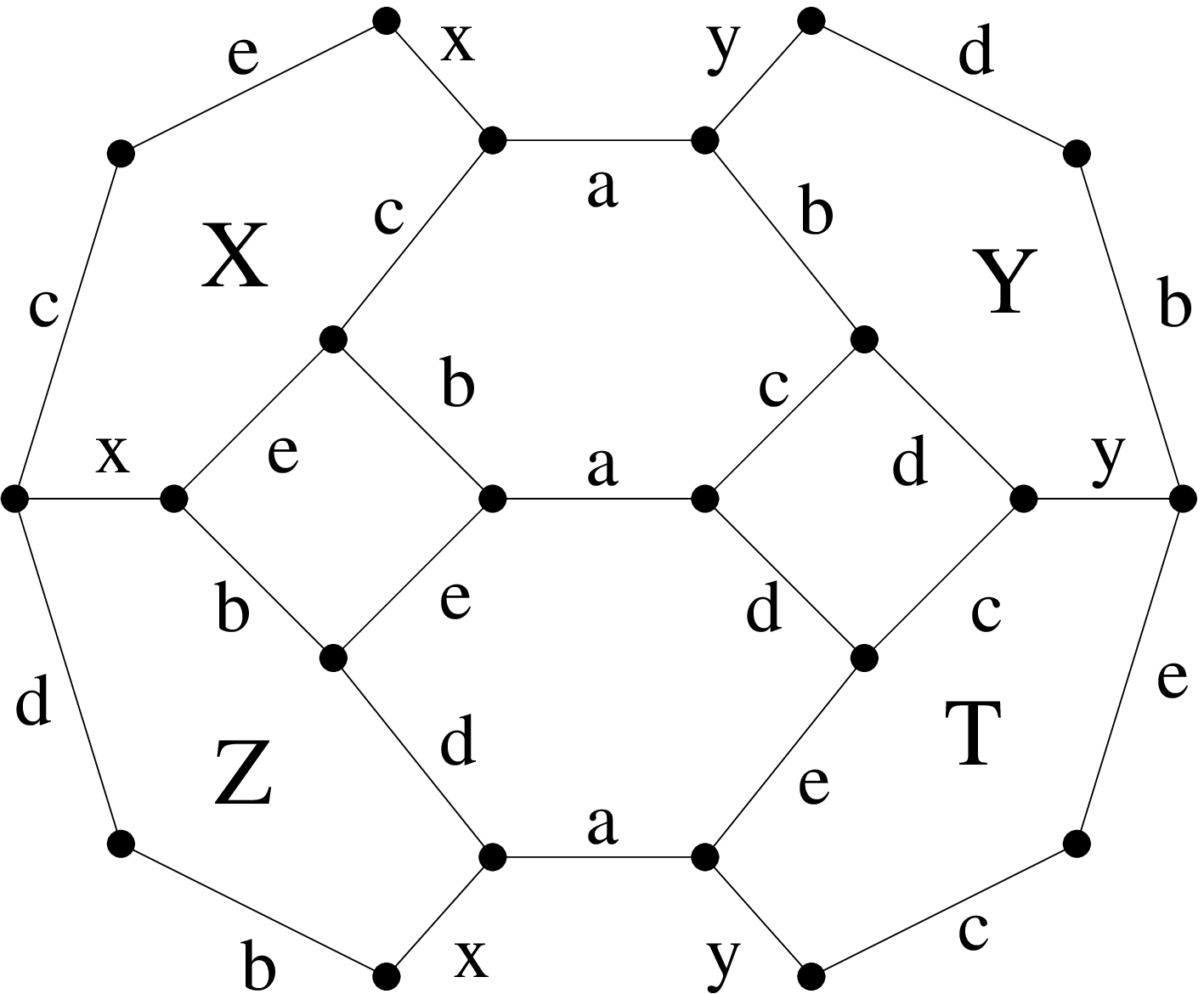,width=5cm}
\end{center}

\begin{lemma}\label{x_and_y_equal}
If in the picture above we have $x=y$ then $\Gm$ is the truncated
Octahedron.
\end{lemma}

\pf Since $x=y$, the two new faces in the center (top and bottom) are
$4$-gonal.
\begin{center}
\epsfig{file=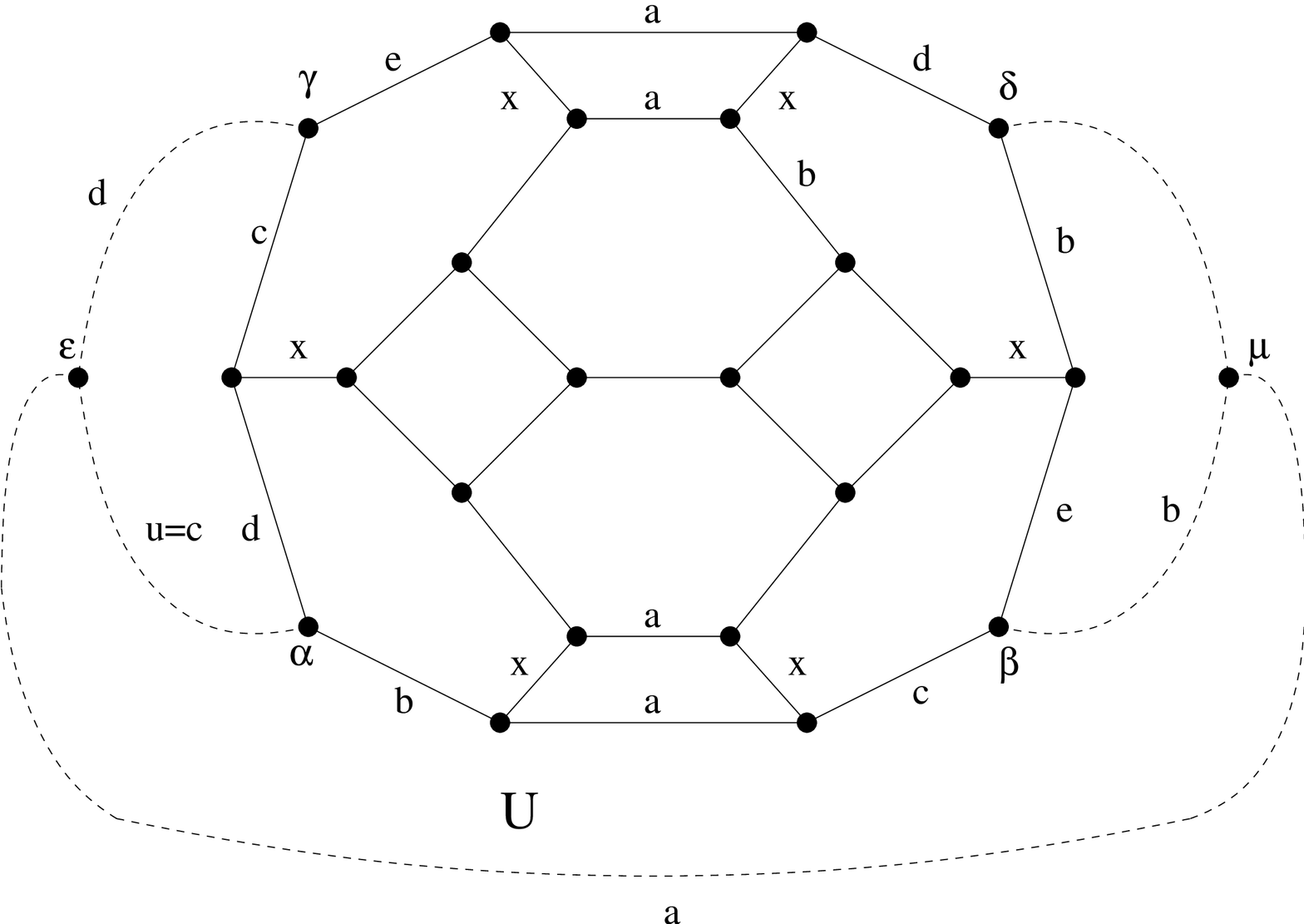,width=7cm}
\end{center}
Consider the unknown label $u$ on the third edge at the vertex $\al$. The
face $U$ has $u$, $b$, $a$, and $c$ as part of its boundary. Hence, $U$ is
$6$-gonal and $u=c$. Similarly, at $\be$, we have an edge labelled $b$, at
$\gm$ an edge labelled $d$, and at $\dl$ an edge labelled $e$. Furthermore,
the new edges at $\al$ and $\gm$ lead to the same vertex $\ep$ and,
similarly, the new edges at $\be$ and $\dl$ lead to a vertex
$\mu$. Finally, we compute that $\ep$ and $\mu$ are connected by an edge
with label $a$.  The resulting graph is trivalent and hence it coincides
with the entire $\Gm$. Manifestly, it is isomorphic to the truncated
Octahedron.\qed

\begin{lemma}\label{x_and_y_different}
Suppose $x\ne y$. Then $\Gm$ is the twisted chamfered Cube.
\end{lemma}

\pf We have the following picture.
\begin{center}
\epsfig{file=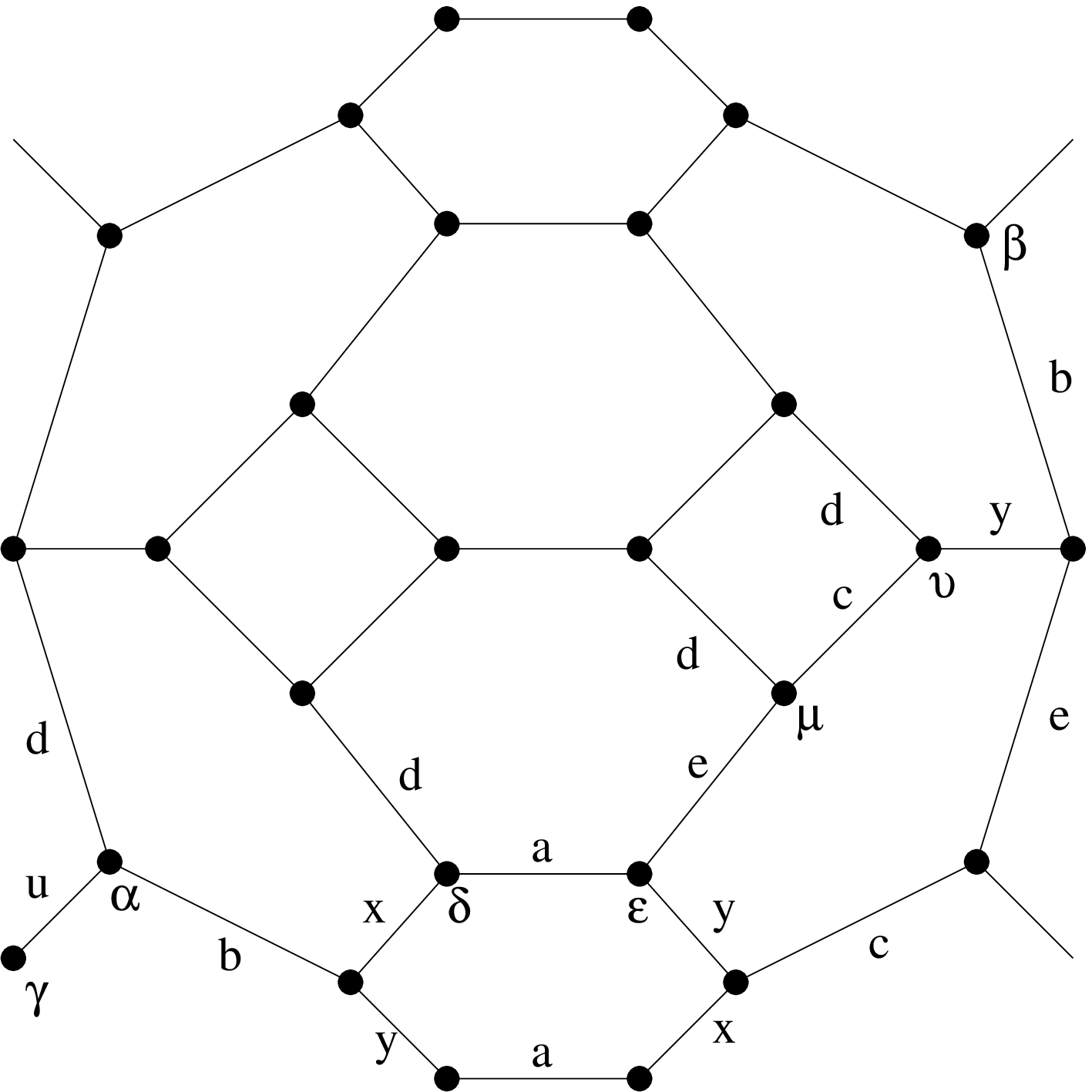,width=5cm}
\end{center}
From this picture we see that the type of every shortest path from $\al$ to
$\be$ is $\{a,c,e,x,y\}$. Thus, $u\in\{a,c,e,x,y\}$. If $u=x$ then $\gm$ is
adjacent to $\dl$ via an edge labelled $b$; a contradiction, since there is
no such edge at $\dl$. So, $u\ne x$. Similarly, if $u=a$ then the shortest
path from $\gm$ and $\ep$ is of type $\{b,x\}$; a contradiction, since none
of these labels appear next to $\ep$. Hence, $u\ne a$. Analogously, looking
at $\gm$ and $\mu$, and at $\gm$ and $\nu$, we establish that $u\ne
e,c$. It follows that $u=y$. Using symmetry, we now know labels on three
more edges, leading to the following picture.
\begin{center}
\epsfig{file=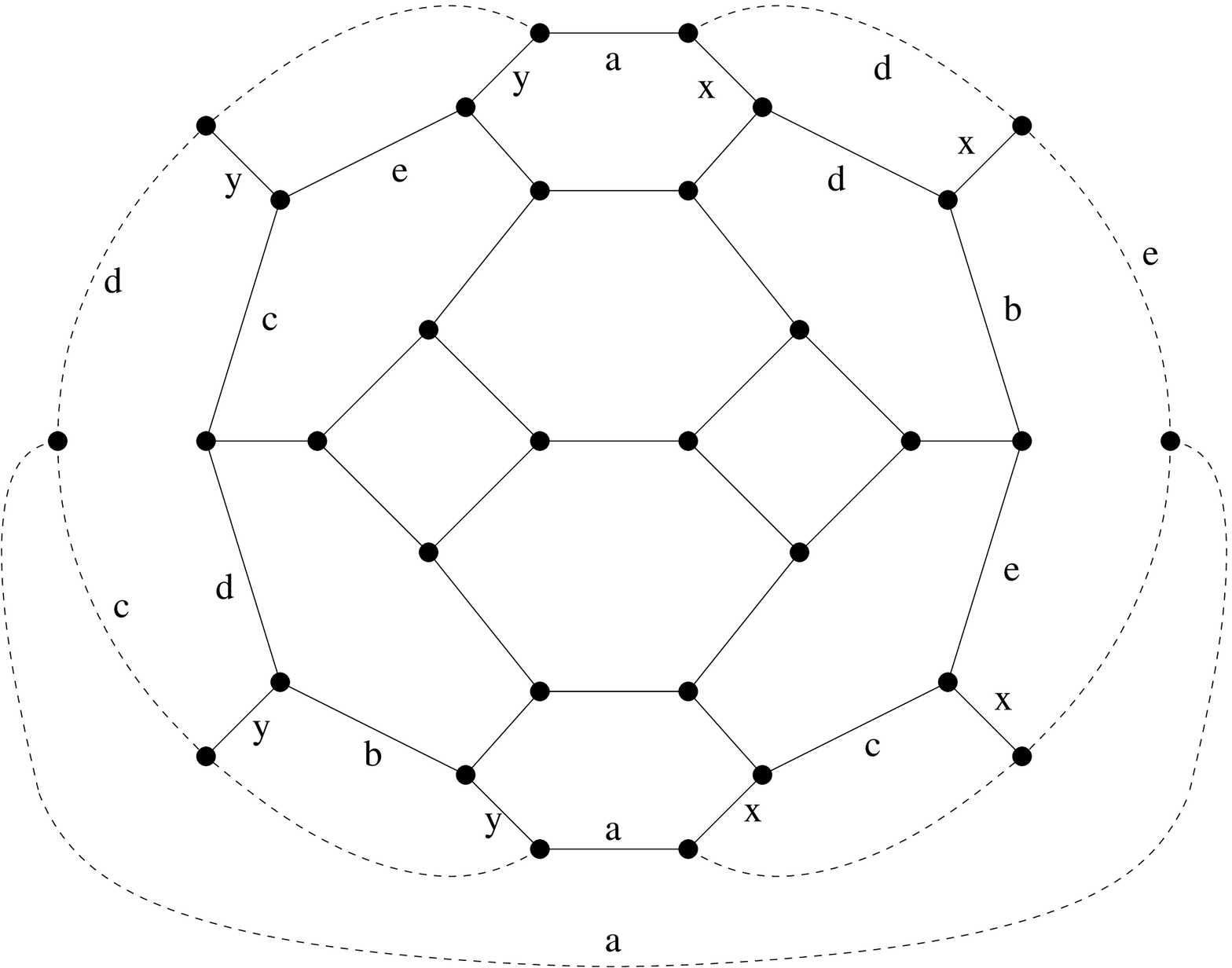,width=8cm}
\end{center}
The consecutive labels $y$, $b$, and $y$ (bottom left) indicate a $4$-gonal
face. Symmetrically, three more $4$-gonal faces appear. Finally, the new
$6$-gonal faces on the left and on the right lead to two more
vertices. These vertices must be adjacent via an edge labelled $a$,
producing a trivalent graph. It follows that this is the entire $\Gm$.
Manifestly, this graph is isomorphic to the twisted chamfered Cube.\qed

\begin{corollary} \label{star}
If $\Gm$ contains the configuration $(\ast)$, but no adjacent $4$-gonal
faces, then $\Gm$ is either the truncated Octahedron, or the twisted
chamfered Cube.\qed
\end{corollary}

We have now looked at all the possibilities arising when $\Gm$ contains the
configuration $(\ast)$. In the remainder of the section {\em we assume
additionally that $(\ast)$ does not occur in $\Gm$}.

We next consider what happens when $\Gm$ contains the following
configuration. 
\begin{center}
\epsfig{file=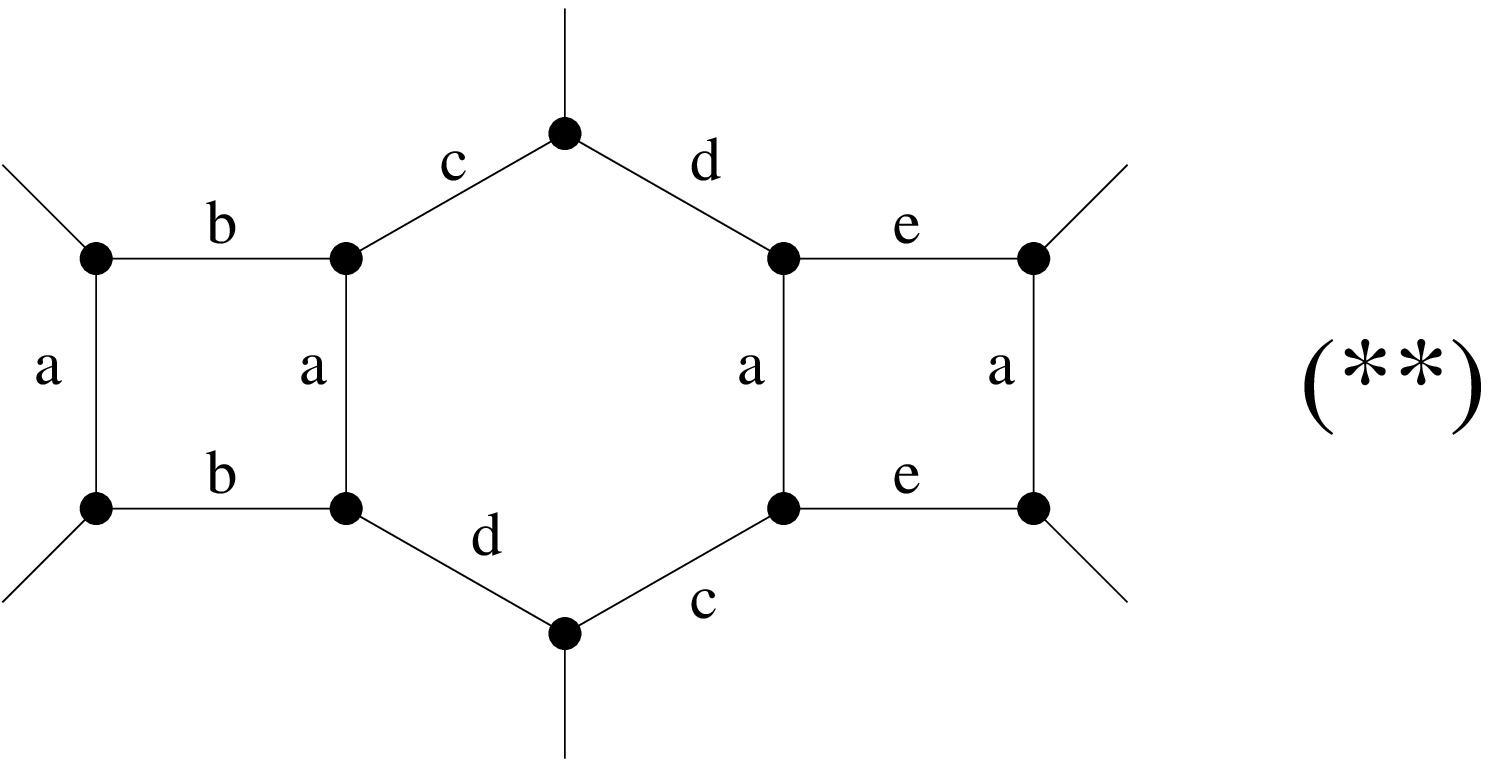,width=5cm}
\end{center}
It follows from Lemma \ref{claw} that the labels $a$, $b$, $c$, $d$, and
$e$ are pairwise distinct.

\begin{lemma}\label{star_star}
Suppose $\Gm$ contains no adjacent $4$-gonal faces and no configuration
$(\ast)$, but it contains $(\ast\ast)$. Then $\Gm$ is isomorphic to the
chamfered Cube.
\end{lemma}

\pf We start with the following picture.
\begin{center}
\epsfig{file=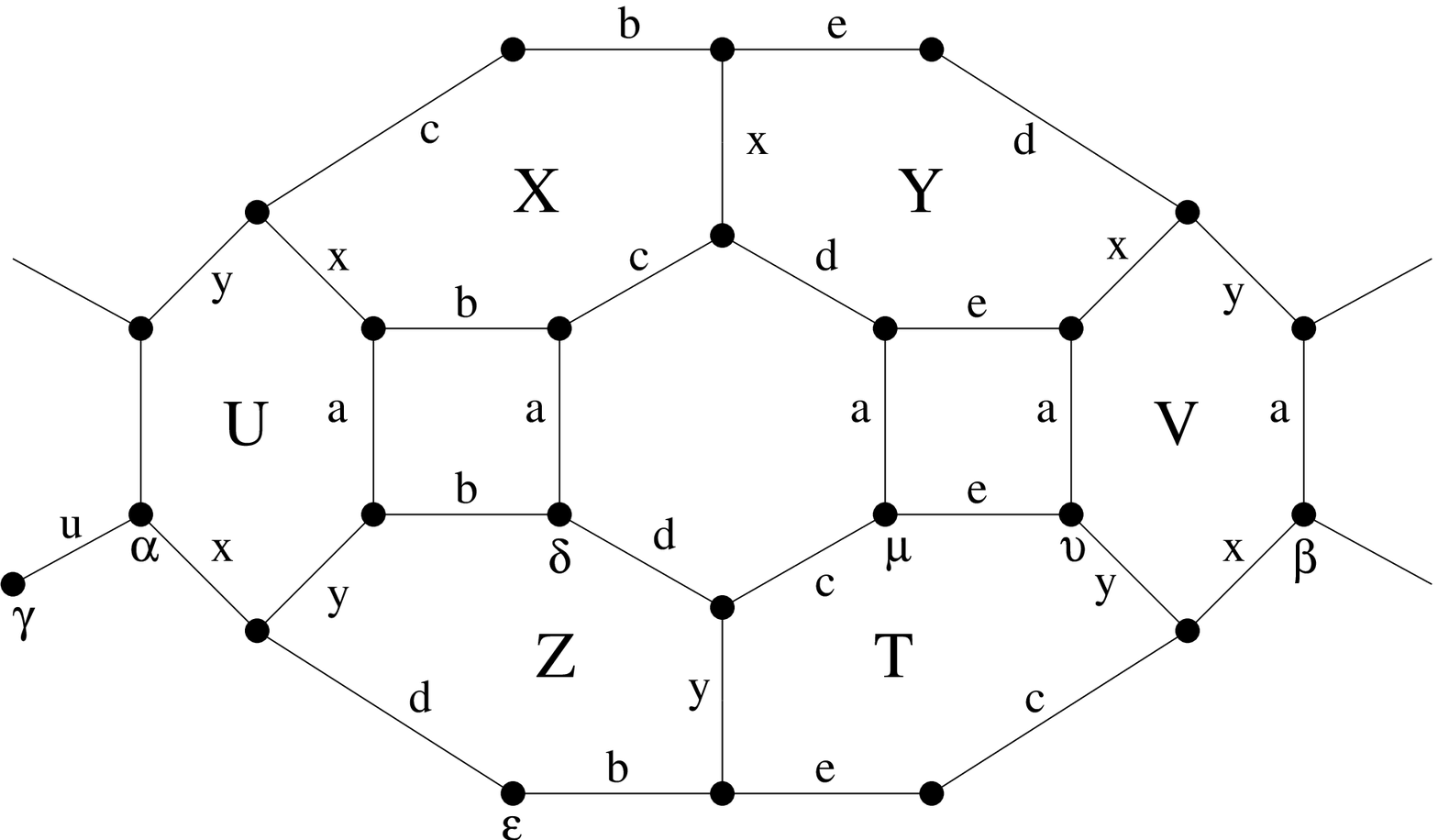,width=7cm}
\end{center}
Observe that the faces $X$, $Y$, $Z$, $T$, $U$, and $V$ must be $6$-gonal
since $\Gm$ contains no adjacent $4$-gonal faces. Thus, only two new
labels, $x$ and $y$, appear, and they cannot be equal to one another, or to
any of the previous labels. Now we compute that the shortest path from
$\al$ to $\be$ must have type $\{b,c,d,e\}$. It follows that the label $u$
on the third edge at $\al$ must be in this set. If $u=b$ then the shortest
path from $\gm$ to $\dl$ has type $\{x,y\}$; a contradiction, since there
no such label next to $\dl$. Hence, $u\ne b$. If $u=d$ then $\gm$ is
adjacent to $\ep$, leading to the configuration $(\ast)$. The contradiction
shows that $u\ne d$. If $u=c$ then the shortest path from $\gm$ to $\mu$
has type $\{b,d,x,y\}$. However, none of these labels appear next to $\mu$.
Therefore, $u=e$. Symmetrically, we determine the labels on three more
edges, ending up with the following picture.
\begin{center}
\epsfig{file=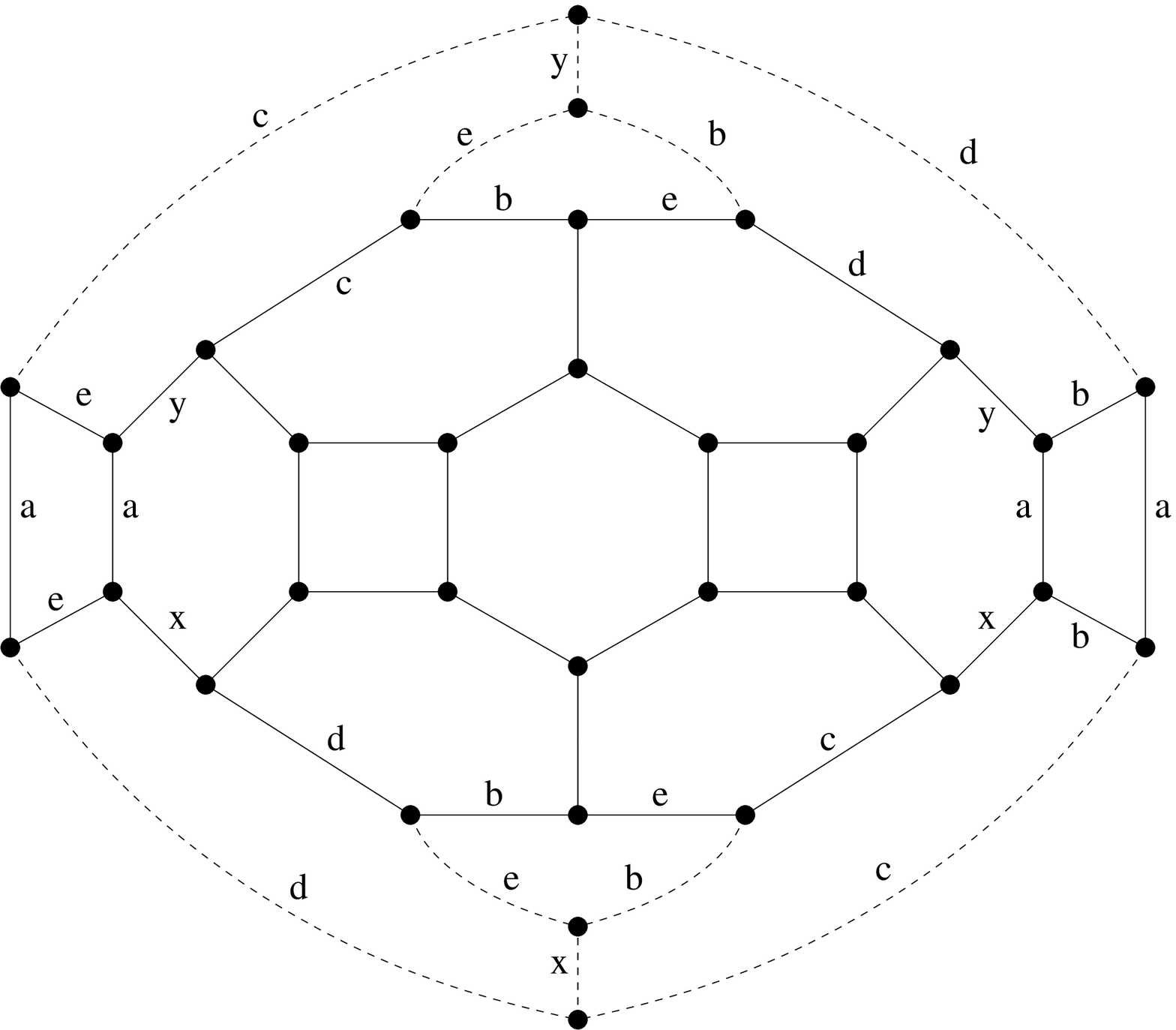,width=9cm}
\end{center}
The consecutive edges labelled $e$, $a$, and $e$ on the left are part of a
$4$-gonal face and, similarly, there is also a $4$-gonal face on the
right. The consecutive edges labelled $e$, $x$, $d$ (bottom left) must be
part of a $6$-gonal face, giving us three more edges labelled $e$, $x$, and
$d$. Symmetrically, there is a $6$-gonal face (bottom right) extending the
sequence of edges $b$, $x$, and $c$. Now the edges labelled $e$ (new), $b$,
$e$, and $b$ (new) must form a $4$-gonal face. Consequently, the two new
edges labelled $x$ are, in fact, one edge. Similar considerations applied
to the top of the picture yield a trivalent graph, which must, therefore,
coincide with $\Gm$. By inspection, $\Gm$ is isomorphic to the chamfered
Cube.\qed

\medskip
At this point we add the assumption that {\em $(\ast\ast)$ also does not
appear in $\Gm$}. Since above we have already met all the five graphs from
the conclusion of Theorem \ref{main}, we will now aim for the final
contradiction.

\begin{lemma} \label{final}
There is no graph $\Gm$ satisfying all of the above assumptions.
\end{lemma}

\pf By contradiction, let $\Gm$ have no adjacent $4$-gonal faces and no
configurations $(\ast)$ and $(\ast\ast)$. We start the analysis of $\Gm$
from the following picture, which must be present in $\Gm$.
\begin{center}
\epsfig{file=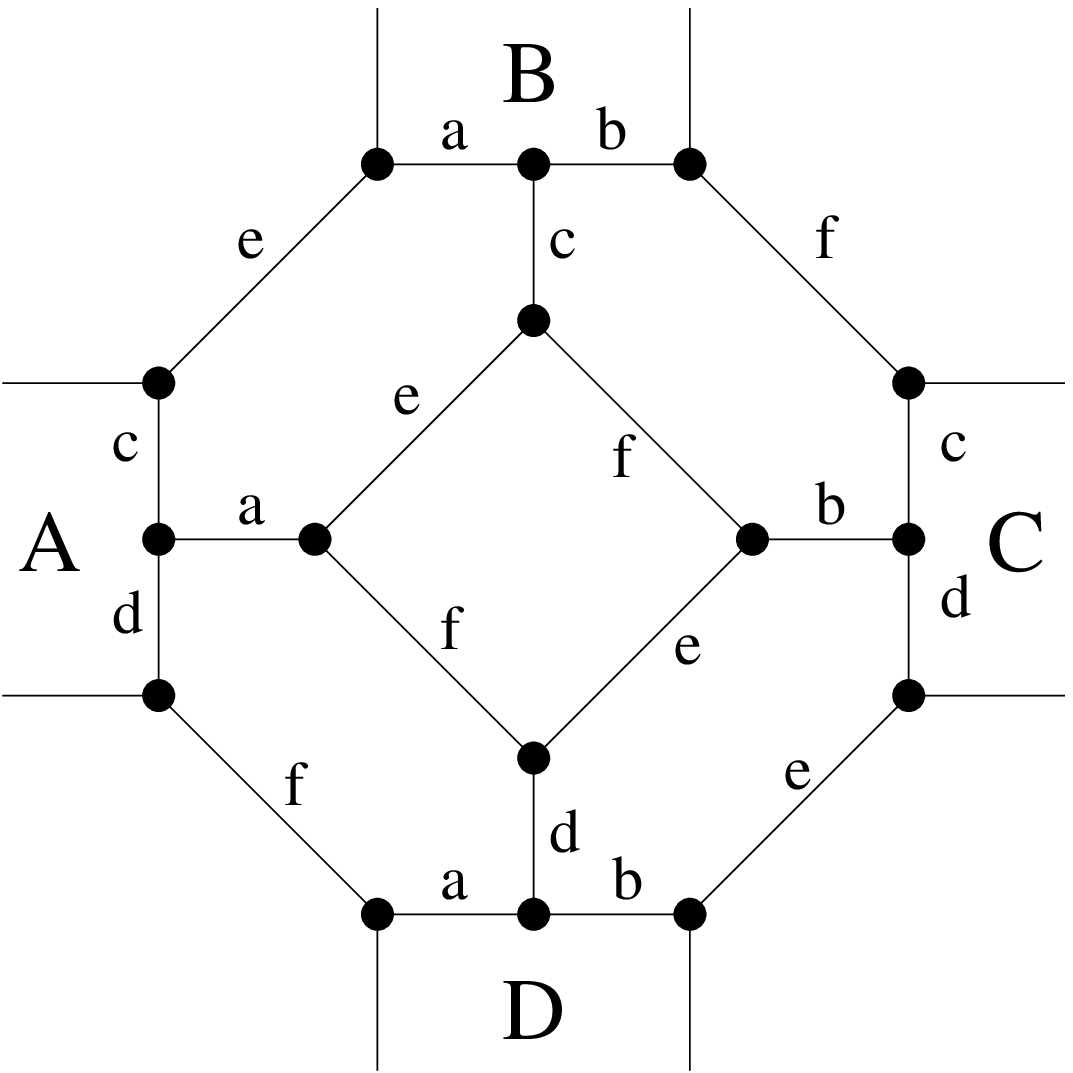,width=5cm}
\end{center}
It is easy to deduce from Lemma \ref{claw} that the labels $a$, $b$, $c$,
$d$, $e$, and $f$ are pairwise distinct. Since $(\ast)$ is not present in
$\Gm$, the faces $A$, $B$, $C$, and $D$ must be $6$-gonal, leading to our
next picture.
\begin{center}
\epsfig{file=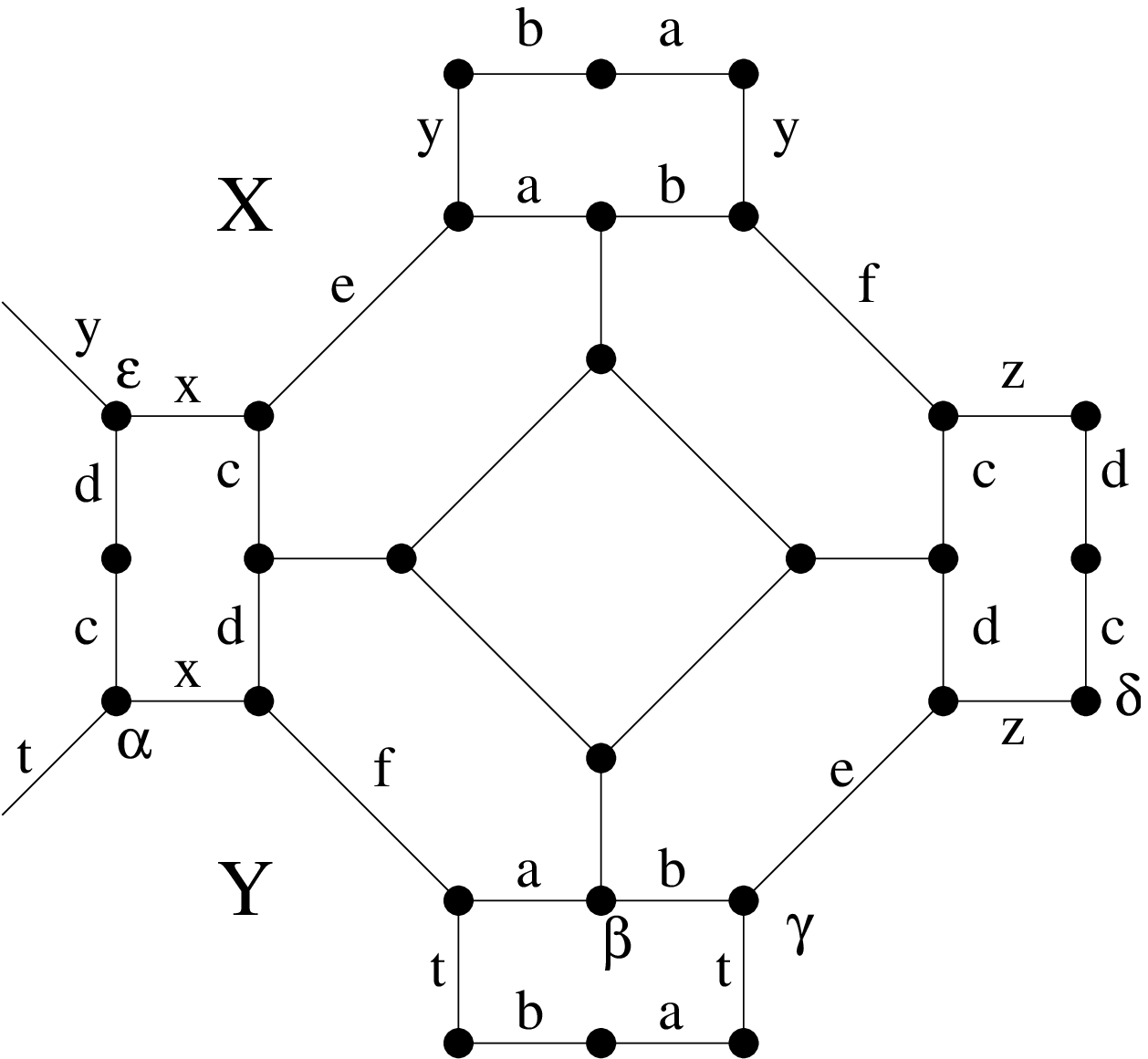,width=7cm}
\end{center}
We claim that the labels $x$, $y$, $z$, and $t$ are new and pairwise
distinct. By symmetry, it suffices to do the check only for $x$. First of
all, $x\not\in\{c,d,e,f\}$ by Lemma \ref{claw}. If $x=a$ then $\al$ must be
adjacent to $\be$ via an edge labelled $f$; clearly, impossible. If $x=b$
the shortest path from $\al$ to $\gm$ has type $\{a,f\}$. However,
$\gm$ does not have these labels next to it (notice that $a\ne t\ne f$ by
Lemma \ref{claw}). Thus, $x\ne b$. The faces $X$ and $Y$ must be $6$-gonal,
since $(\ast\ast)$ does not occur in $\Gm$. In particular, this means that
$x\ne y,t$. Furthermore, this also means that the third edge incident to
$\al$ carries the label $t$ and the third edge at $\ep$ carries the label
$y$. We can now, finally, establish that $x\ne z$. Indeed, if $x=z$ then
the shortest path from $\al$ to $\dl$ must have type $\{a,b,e,f\}$, in
contradiction with the fact that none of the labels next to $\al$ belong to
this list. Thus, indeed, no two labels in $\{a,b,c,d,e,f,x,y,z,t\}$ are
equal.

It is time for our final picture.
\begin{center}
\epsfig{file=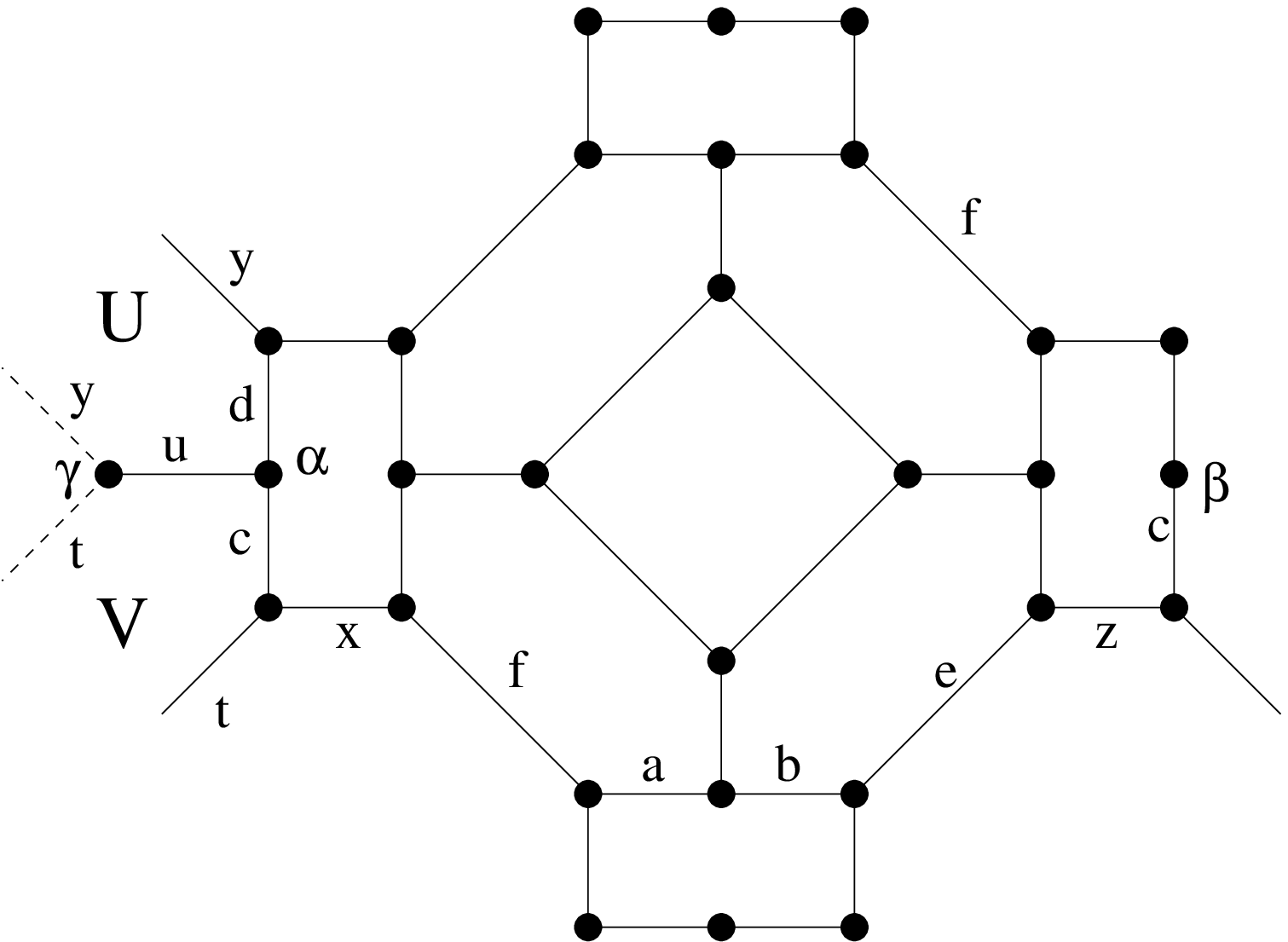,width=8cm}
\end{center}
Let $u$ be the label on the third edge at the vertex $\al$ from this
picture. Since the shortest path from $\al$ to $\be$ must have type
$\{a,b,e,f,x,z\}$, that path must go via the edge labelled $u$, and $u$
must be on the above list. In particular, $t\ne u\ne y$. This means that
the faces $U$ and $V$ are $6$-gonal. This, in turn, implies that the
remaining two edges at the vertex $\gm$ have labels $t$ and $y$. This gives
the final contradiction: Since neither $t$, nor $y$ belong to
$\{a,b,e,f,x,z\}$, the shortest path from $\al$ to $\be$ cannot pass
through $\gm$, either.
%\qed

\medskip
This concludes the proof of Theorem \ref{main}. \qed

\section{Remarks and possible extensions}

A {\em $t$-embedding} of a graph $G$ is a mapping 
$\phi:G\mapsto H_m$, such that it holds:
\begin{equation*}
d(\phi(x),\phi(y))=d_G(x,y)\mbox{~~~if~~~}d_G(x,y)\leq t.
\end{equation*}

\begin{remark}
(i) Clearly, the five graphs $4_n$ from Theorem \ref{main} have the following 
isometric embedding:
Cube into $H_3$, $6$-gonal prism into $H_4$,
truncated Octahedron into $H_6$;
chamfered Cube and its twist into $H_7$.
Only the first four are zonohedra. Only the first three are 
space-fillers (moreover, Voronoi polyhedra of lattices $Z_3$,
$A_2 \times Z_1$ and $A_3^*$, respectively).

(ii) Any dual $4_n$ with $n>4$ is not $l_1$-embeddable (moreover,
not $3$-embeddable), because it contains, always isometric,
the non $5$-gonal subgraph (of diameter $3$), depicted on Figure 13.3
(right hand side) of \cite{DGS}.
\end{remark}

\begin{remark}\label{REM2}
(i) No $3_{n}$, except of Tetrahedron, is $l_1$-embeddable;
moreover, it is not $4$-embeddable.
%Unique $3_8$ and unique $3_{12}$ are $2$- but not $3$-embeddable;
%any $3_n$ with $n>12$ is not $4$-embeddable.

In fact, the self-dual $3_4$ (Tetrahedron) have isometrical
embeddings into $\frac{1}{2}H_3$ and $\frac{1}{2}H_4$.
The unique case of $3_n, n \ge 4,$, having non-isolated
triangles, is the unique (and it is non $3$-connected one) $3_8$.
Easy to check that
this graph of diameter $3$ is not $5$-gonal; so, it is not $l_1$-embeddable.
%(Apropos, it is $2$-embeddable and its dual embeds isometrically into
%$\frac{1}{2}H_6$.)
The unique polyhedron $3_{12}$ (truncated Tetrahedron) also has diameter
$3$ and is not $5$-gonal.

Now, the non $5$-gonal configuration (of diameter $4$), depicted on Figure 13.3
(left hand side) of \cite{DGS}, show that any $3_n$, having it as an
{\em isometric} subgraph, is not $5$-gonal.

For $n>12$, the only $3_n$, for which the non $5$-gonal subgraph from
Figure 13.3 of \cite{DGS} is not isometric, are those of class (ii) from
Theorem 5.1 of \cite{zig2}, classifying all $3_n$.
Other non $5$-gonal subgraph, of diameter $4$ again, will be isometric
for all those $3_n$.

(ii) All known $l_1$-embeddable $(3_n)^*$ embeds into $\frac{1}{2}H_m$
and have (see \cite{DG2})\\
$(n,Aut,m)$=$(4, T_d; 3)$, $(4, T_d; 4)$, $(8, D_{2h}; 6)$,
$(12, T_d; 7)$, $(16, D_{2h}; 8)$, $(16, T_d; 8)$, $(28, T; 10)$,
$(36, T_d; 11)$.
\end{remark}

\begin{remark}(see \cite{DG2})\label{REM3}
(i) All known $l_1$-embeddable $5_n$  embeds into $\frac{1}{2}H_m$ and have\\
$(n, Aut; m)$=$(20, I_h; 10)$, $(26, D_{3h};12), (40, T_d; 15)$, $(44, T;16)$,
$(80, I_h; 22)$.

(ii) All known $l_1$-embeddable $(5_n)^*$ embeds into $\frac{1}{2}H_m$ and have\\
$(n, Aut; m)$=$(12, I_h; 6)$, $(28, T_d; 7)$, $(36, D_{6h}; 8)$, $(60, I_h; 10)$.
\end{remark}

{\bf Conjecture}

All $l_1$-embeddable graphs $3_n$, $5_n$ and their duals are those given
in the above Remarks \ref{REM2} and \ref{REM3}.

A {\em zone} in a bipartite plane graph is a circuit 
$(F_i)_{1\leq i\leq 
h}$ of faces, such that each face $F_i$ is adjacent to $F_{i-1}$ and 
$F_{i+1}$ (in cyclic order) on opposite edges.

\begin{proposition}
If a bipartite plane graph is $t$-embeddable, has faces of gonality
at most $t$ only and such that a shortest path between two vertices
of a face $F$ is included in $F$, then none of its zones is
self-intersecting.
\end{proposition}
\pf From the $t$-embedding condition and the face condition, one obtains 
that the labels corresponding to opposite edges are identical.
Also, two identical labels can occur on a face only on opposite edges.
So, all zones are not self-intersecting. \qed

Note that it is proved in \cite{DGS}, p.~151 that a bipartite plane graph
without self-intersecting zones is {\em strictly admissible},
i.e. admits an $1$-embedding in $H_m$, such that faces are mapped to isometric
cycles of $H_m$.

The above Proposition provides an efficient tool for selecting the 
possible bipartite plane graphs, which are $t$-embeddable with $2t$
being the maximal gonality of the faces.
A computation for graphs $4_n$ up to $n=210$ vertices gave:
\begin{enumerate}
\item all 5 graphs in the Main Theorem and
\item all graphs (up to $210$ vertices) with $Aut=O_{h}$.
\end{enumerate}
All $4_n$ with $Aut=O_h$ are exactly those, which come from Cube
by {\em Goldberg-Coxeter construction} (see \cite{DuDe} for definitions 
and details); one can see that none of them have a self-intersecting 
zone.

We expect that those graphs and the $5$ graphs of Main Theorem are only
$3$-embeddable graphs $4_n$.

The only known $l_1$-embeddable $3$-valent plane graphs 
with only $4$- and $2m$-gonal faces, $m>3$, are Cube and $2m$-gonal 
prism (an isometric subgraph of $H_m$).
%But we expect many such 
%$3$-embeddable graphs; for example, the number of such graphs with 
%$m=4$, at most $56$ vertices and without self-intersecting railroads, is 
%$19$.

\end{document}